\documentclass[a4paper, twoside, 10pt]{article}
\usepackage{eurosym, xcolor}
\usepackage{amsmath,amsthm,amsfonts,latexsym,amscd,amssymb, enumerate}
\usepackage{hyperref}
\numberwithin{equation}{section}
\input xypic
\newtheorem{theorem}{Theorem}[section]
\newtheorem{lemma}{Lemma}[section]
\newtheorem{corollary}{Corollary}[section]
\newtheorem{proposition}{Proposition}[section]
\newtheorem{definition}{Definition}[section]
\newtheorem{example}{Example}[section]
\theoremstyle{remark}

\date{}
\DeclareMathOperator{\R}{\mathbb{R}}%
\title{\textbf{Clairaut anti-invariant Riemannian maps with K\"ahler and Ricci soliton structures}}
\author{Jyoti Yadav and Gauree Shanker\thanks{corresponding author, Email: gauree.shanker@cup.edu.in}}

\begin{document}
	\maketitle
	\begin{abstract}
		The aim of this article is to explore the Clairaut anti-invariant Riemannian maps from/to K\"ahler manifolds admitting Ricci solitons. We find the curvature relations and calculate the Ricci tensor under different conditions. We discuss the condition under which range space becomes $\alpha$-Ricci soliton. We obtain conditions for the range and kernel spaces of these maps to be Einstein. Next, we find the scalar curvature for range space. Further, we give the relation between Ricci curvature and Lie derivative under these maps. Moreover, we find the condition for a vertical potential vector field on target manifold to be  conformal vector field on range space of these maps. Finally, we give non-trivial examples of such maps.
	\end{abstract}
	\noindent\textbf {M. S. C. 2020:} 53B20, 53B35.
	
	\section{Introduction}\label{sec1}
	The submersions and immersions play an important role in the study of smooth maps between smooth manifolds. In the theory of Riemannian manifolds, Riemannian maps were introduced by Fischer \cite{S27} in 1992 as a generalization of isometric immersion and Riemannian submersion. The notable property of Riemannian maps is that these maps satisfy the generalized eikonal equation, which connects the geometric and physical optics.
	
	The geometry of Riemannian submersions was discussed by Falcitelli et al.  \cite{S3}. As a generalization of holomorphic Riemannian submersions, \c{S}ahin \cite{S31} introduced holomorphic Riemannian maps. As a generalization of anti-invariant Riemannian submersions, \c{S}ahin \cite{S32} introduced anti-invariant Riemannian maps from almost Hermitian manifolds to Riemannian manifolds.
	
	The idea of Clairaut's Theorem is based on geodesics on the surface of revolution. This theorem states that for any geodesic $\gamma:I\subset \mathbb{R}\rightarrow S$ on a surface of revolution, $r\sin\theta$ is constant along $\gamma$, where $\theta(s)$ is the angle between $\gamma(s)$ and the meridian curve through $\gamma(s)$, $s \in I$ and  $r : S \rightarrow \mathbb{R}^+$ is the distance of a point of $S$ from the axis of rotation. The authors of \cite{N2, Meena-Zawadzki} generalized this idea on submersion theory and defined Clairaut submersions. Further, this idea has been generalised to Riemannian maps as necessary and sufficient conditions for Riemannian maps to be a Clairaut Riemannian maps obtained in \cite{S2, S35, KSS}. Recently, Clairaut anti-invariant Riemannian maps from/to K\"ahler manifolds have been studied in \cite{S35, S20}.
	
	Ricci solitons are the well-known object in Riemannian geometry. It comes out as the solution of Einstein's field equations within the framework of general relativity. Since the equation of Ricci soliton has relations with string theory, it has also become an attractive subject of study for physicists. \cite{S33} Ricci soliton appears during the analysis of the Ricci flow \\
	$$\frac{\partial g(t)}{\partial t} = -2g(t),$$
	where  $g(t)$ is a one-parameter family of metrics on a  manifold.
	
	A smooth vector field $X'$ on a Riemannian manifold $(M, g_M)$ is called a conformal vector field if the Lie derivative of the metric $g_M$ with respect to $X'$ can be expressed as a multiple of a function to the metric $g_M$, and this function is called the conformal factor of the conformal vector field $X'$ \cite{N6}.

	A connected Riemannian manifold $(M, g_M)$ is called $\alpha$ - Ricci soliton \cite{N5}, denoted by $(M, g_M, \xi, \lambda)$, if it satisfies the following equation
	\begin{equation}\label{Ricci eqn}
		\frac{1}{2}(L_{\xi}g_M)(X, Y) + \alpha Ric(X, Y) +\lambda g_M(X, Y)= 0,
	\end{equation}
	where $\alpha, \lambda\in\R,$ $\xi$ is a potential vector field, $L_{\xi}g_M$ denotes the Lie derivative of metric tensor $g_M$ with respect to $\xi$, Ric denotes the Ricci curvature.
	The $\alpha$ - Ricci soliton is called shrinking, steady or expanding according to $\lambda < 0, \lambda = 0$ or $\lambda > 0,$ respectively. If $\alpha = 1,$ then \eqref{Ricci eqn} reduces to Ricci soliton.
	If $\xi$ is gradient of some function $f$ which is defined on $M$, i.e., $\xi = grad f$, $f\in C^{\infty}(M)$, then from \eqref{Ricci eqn}, we get\\
	\begin{equation*}
		Hess^f(X, Y)+ Ric(X, Y) + \lambda g_M(X, Y) = 0,
	\end{equation*}
	where $Hess^f$ is hessian form. Then $(M, g_M, \nabla f, \lambda)$ is called gradient Ricci soliton and the function f is called potential function. 
	Pigola et al. \cite{S34} gave the idea of almost Ricci soliton by using $\lambda$ as a variable function rather than constant function. In this case $(M, g_M, \xi, \lambda)$ is called an almost Ricci soliton.
	
	The authors of \cite{S39, Meric_2020, Meena-Conformal} study about submersions whose total manifolds admit a Ricci soliton. Yadav, Meena \cite{N7, S51} and Gupta et al. \cite{GGupta_2022} have studied Riemannian maps whose total (base) manifold admit a Ricci soliton. In \cite{S35, S36, KSS}, Meena and Yadav have studied Clairaut Riemannian maps whose total (base) manifolds admit a Ricci soliton. Watson \cite{S6} introduced an almost Hermitian submersion and  G\"und\"uzalp \cite{Gunduzalp_2020} studied an almost Hermitian submersion \cite{S52} whose total manifold admits Ricci soliton.

	In this paper, we explore Clairaut anti-invariant Riemannian maps from/to K\"ahler manifolds admitting Ricci solitons. The paper is divided into five sections. In section 2, we recall some definitions and results which are useful to current study. In section 3, we study Clairaut anti-invariant Riemannian maps from a K\"ahler manifold admitting Ricci soliton to a Riemannian manifold and give non trivial examples for the existence. In section 4, we study Clairaut anti-invariant Riemannian maps from a Riemannian manifold to a K\"ahler manifold admitting Ricci soliton and construct some non-trivial examples. In section 5, we give scope for further studies based on current study.
	
	
	\section{Preliminaries}
	Let $F$ be a smooth map between two Riemannian manifolds $(M, g_M)$ and $(N, g_N)$ of dimension $m, n,$ respectively. Let $\mathcal{V}_p = kerF_{*p}$ at $p\in M,$ denotes vertical distribution or kernel space of $F_*$ and $\mathcal{H}_p = (kerF_{*p} )^\perp$ in $T_pM$ is the orthogonal complementary space of $\mathcal{V}_p$. 
	Then, we decompose the tangent space $T_pM$ at $p\in M$ in the following manner.\\
	$$ T_pM = (kerF_{*p}) \oplus (kerF_{*p} )^\perp = \mathcal{V}_p \oplus \mathcal{H}_p.$$
	Let the range of $F_*$ be denoted by $rangeF_{*p}$ at $p\in M$, and $(rangeF_{*p})^\perp$ be the orthogonal complementary space of $rangeF_{*p}$ in the tangent space $T_{F(p)}N$ of $N$ at $F(p)\in N$. Since $rankF <\min\{m,n\}$, we get $(rangeF_{*p})^\perp\neq \{0\}$. Thus, the tangent space $T_{F(p)}N$ of $N$ at $F(p)\in N$ has the following decomposition:
	\begin{equation*}
		T_{F(p)}N = (rangeF_{*p})\oplus (rangeF_{*p})^\perp.
	\end{equation*}
	A smooth map $F$ is said to be a Riemannian map at $p\in M$ if the horizontal restriction
	$F^h_{p}:(kerF_{*p})^\perp\rightarrow (rangeF_{*p})$ is a linear isometry between the inner product spaces $((kerF_{*p})^\perp, g_M{(p)}| (kerF_{*p})^\perp)$ and $((rangeF_{*p}, g_N(q)|(rangeF_{*p})), q = F(p)$. In other words, $F_*$ satisfies the equation \cite{S30}\\
	\begin{equation}{\label{N1}}
		g_N(F_{*}X, F_{*}Y) = g_M(X, Y)~~ \forall X, Y \in\Gamma(kerF_*)^\perp.
	\end{equation}
	For Riemannian submersion, O'Neill \cite{S10} defined fundamental tensor fields $T$
	and $A$ by
	\begin{equation}\label{neill A}
		A_{X'}{Y'} = \mathcal{H}\nabla^M_{\mathcal{H}X'}\mathcal{V}Y' + \mathcal{V}\nabla^M_{\mathcal{H}X'}\mathcal{H}Y',\\
	\end{equation}
	\begin{equation}\label{neill T}
		T_{X'}{Y'} =  \mathcal{H}\nabla^M_{\mathcal{V}X'}\mathcal{V}Y' +\mathcal{V}\nabla^M_{\mathcal{V}X'}\mathcal{H}Y',\\
	\end{equation}
	where $X', Y'\in\Gamma(TM), \nabla^M$ is the Levi-Civita connection of $g_M$ and $\mathcal{V}, \mathcal{H}$ denote the projections to vertical subbundle and horizontal subbundle,
	respectively. For any $X'\in \Gamma(TM)$, $T_{X'}$ and $A_{X'}$ are skewsymmetric
	operators on $(\Gamma(TM), g_M)$ reversing the horizontal and the vertical distributions.
	One can check easily that $T$ is vertical, i.e., $T_{X'} = T_{\mathcal{V} X'}$, and $A$ is horizontal, i.e., $A = A_{\mathcal{H}X'}$.\\
	We observe that the tensor field $T$ is symmetric for vertical vector fields, 
	and tensor field $A$ is anti-symmetric for horizontal vector fields. 
	Using  equation \eqref{neill A} and \eqref{neill T}, we have the following Lemma.\\
	\begin{lemma} \cite{S10} Let  $X, Y\in\Gamma (kerF_{*})^\perp$ and $V,W \in \Gamma(kerF_{*}).$ Then
		\begin{equation*}
			\nabla_{V}W = T_{V}W + \hat{\nabla}_{V}W,
		\end{equation*}
		\begin{equation*}
			\nabla_{V}X = \mathcal{H}\nabla_{V}X + T_{V}X,
		\end{equation*}
		\begin{equation*}
			\nabla_{X}V = A_{X}V +\mathcal{V}\nabla_{X}V,
		\end{equation*}
		\begin{equation*}
			\nabla_{X}Y = \mathcal{H}\nabla_{X}Y + A_{X}Y.
		\end{equation*}
		If $X$ is a basic vector field. then $\mathcal{H}\nabla_{V}X = A_X V.$
	\end{lemma}
	For any vector field $X'$ on $M$ and any section $D$ of $(rangeF_*)^\perp$, we define the connection $\nabla^{F\perp}$ as the orthogonal projection of $\nabla_{X'} D$ on $(rangeF_*)^\perp.$ One can easily see that $\nabla^{F\perp}$ is a linear connection on $(rangeF_*)^\perp$ with $\nabla^{F\perp}g_N = 0.$
	For a Riemannian map $F$, we have  \cite{S23}
	\begin{equation}\label{S_V}
		\nabla^{N}_{F_{*}X}D = -S_D{F_{*}X} + \nabla^{F \perp}_ {X} D,
	\end{equation}
	where $S_D{F_{*}X}$ is the tangential component and  $\nabla^{F \perp}_ {X} D$ is the orthogonal component  of $\nabla^{N}_{F_{*}X}D$.
	It can be easily seen that $\nabla^{N}_{F_{*}X}D$ is obtained from the pullback connection of $\nabla^{N}$. Thus, at $p\in M$, we have 
	$\nabla^{N}_{F_{*}X}D(p) \in T_{F(p)}N, S_D{F_{*}X}\in F_{*p}(T_{p}M)$ and $\nabla^{F \perp}_ {X} D\in (F_{*p}(T_{p}M))^\perp$. It follows that $S_D{F_{*}X}$ is bilinear in $D$,  $F_{*}X$ and $S_D{F_{*}X}$  at $p$ depends only on $D_{p}$ and $F_{*p}X_{p}$.\\  By direct computations, we obtain 
	\begin{equation*}\label{Sv}
		g_N(S_D{F_{*}X, {F_{*}Y}}) = g_N \big(D, (\nabla F_{*})(X, Y)\big)\\ 
	\end{equation*} 
	for all X, Y $\in \Gamma(kerF_{*})^\perp~and ~V\in\Gamma(rangeF_{*})^\perp.$
	Since $\nabla F_{*}$ is symmetric, it follows
	that $S_D$ is a symmetric linear transformation of $range F_{*}$.  \\\
	Let $F : (M, g_M) \rightarrow (N, g_N)$ be a smooth map between smooth manifolds. The second fundamental form of $F$ is the map \cite{S37}
	\begin{equation*}
		\nabla F_{*} : \Gamma(TM)\times \Gamma(TM) \rightarrow \Gamma_{F}(TN)\\
	\end{equation*}
	defined by
	\begin{equation*}\label{SFF}
		(\nabla F_{*}) (X', Y') = \nabla^{N_{F}}_{X'}{F_{*}}Y' - F_{*}(\nabla^M_{X'} Y'),
	\end{equation*}
	where  $\nabla^M$ is a linear connection on $M$.\\
		\begin{lemma}\label{Lemma}\cite{S38} Let $F : (M, g_M) \rightarrow (N, g_N)$ be a Riemannian map between Riemannian manifolds. Then $F$ is an umbilical Riemannian map if and only if $$(\nabla F_*)(X, Y ) = g_M(X, Y )H'$$ for all $X, Y \in \Gamma (kerF_*)^\perp$ and $H'$ is nowhere zero mean curvature vector field on $(rangeF_*)^\perp$. 
	\end{lemma}
	In \cite{S1}, it is proved that an almost Hermitian manifold $(M, g_M)$ admits a
tensor $J$ of type (1,1) such that $J^2 = -I$ and
\begin{equation}\label{COHM}
	g_M(JX', JY') = g_M(X', Y')
\end{equation}
for all $X', Y' \in \Gamma(TM).$ An almost Hermitian manifold $M$ is called a K\"ahler manifold if
\begin{equation}\label{COKM}
	(\nabla^M_{X'}J)Y' = 0\\
\end{equation}
for all $X', Y' \in \Gamma(TM)$, where $\nabla^M$ is Levi-Civita connection on $M$.
\subsection{Some results on Clairaut anti-invariant Riemannian maps from K\"ahler manifolds}
In this subsection, we recall some definitions and results which will be used in section 3. The first definition deals with Clairaut Riemannian map for source manifold which was defined by \c{S}ahin. 
\begin{definition}
\cite{S2} A Riemannian map $F:(M, g_M)\rightarrow (N, g_N)$ between Riemannian manifolds is called a Clairaut Riemannian map if there is a function
$\tilde{r} : M \rightarrow \mathbb{R}^+$ such that for every geodesic, making angles $\theta$ with the horizontal subspaces, $\tilde{r} \sin \theta$ is constant.	
\end{definition}
\begin{theorem}\cite{S2}\label{CCRmp}
	Let $F : (M, g_M) \rightarrow (N, g_N)$ be a Riemannian map between Riemannian manifolds with connected fibers, then $F$ is a Clairaut Riemannian map with $\tilde{r} = e^f$ if and only if each fiber is umbilical and has mean curvature vector field $H = -gradf$, where $gradf$ is the gradient of the function $f$ with respect to $g_M$.	
\end{theorem}
	\begin{definition}\cite{S32}
	Let $F :(M, g_M, J)\rightarrow (N, g_N)$ be a Riemannian map from an almost Hermitian manifold  to a Riemannian manifold. We say that $F$ is an anti-invariant Riemannian map if 
	$$J(ker F_*) \subset (ker F_*)^\perp.$$
	We denote the orthogonal complementary subbundle to $J(ker F_*)$ in $(ker F_*)^\perp$ by $\mu$.	
\end{definition}
Next, let $F: (M, g_M, J)\rightarrow (N, g_N)$ be an anti-invariant Riemannian map from a K\"ahler manifold to a Riemannian manifold. Then, for any $X \in\Gamma(ker F_*)^\perp,$ we have
\begin{equation}\label{JX}
	JX = BX + CX,	
\end{equation}
where $BX\in \Gamma(ker F_*)$ and $CX \in \Gamma(\mu).$ Note that if $\mu = 0,$ then $F$ is called a Lagrangian Riemannian map.\\
Now, we recall the definition of Clairaut anti-invariant Riemannian map for source manifold.
\begin{definition}
\cite{S20} An anti-invariant Riemannian map from K\"ahler manifold to
a Riemannian manifold is called Clairaut anti-invariant Riemannian map if
it satisfies the definition of Clairaut Riemannian map for source manifold.
\end{definition}

\subsection{Some results on Clairaut anti-invariant Riemannian maps to K\"ahler manifolds}
In this subsection, we recall some definitions and results which will be used in section 4. The first definition is about Clairaut Riemannian map for target manifold.
	\begin{definition}\label{CRMps}\cite{S35}
	A Riemannian map $F : (M, g_M) \rightarrow (N, g_N)$ between Riemannian manifolds is called a Clairaut Riemannian map if there is a function $\tilde{s}: N \rightarrow \mathbb{R}^+$ such that for every geodesic $\sigma$ on $N,$ the function 
	$(\tilde{s} o \sigma) sin \omega(t)$  is constant, where $F_*X \in \Gamma(rangeF_*)$ and $D \in \Gamma(range F_*)^\perp$ are the vertical and horizontal components of $\dot\sigma(t)$, and $\omega(t)$ is the angle between $\dot\sigma(t)$ and $D$ for all t.
\end{definition}
	\begin{theorem}\cite{S35}\label{NSC}
		Let $F : (M, g_M ) \rightarrow (N, g_N)$ be a Riemannian map between Riemannian manifolds such that $(rangeF_*)^\perp$ is totally geodesic and $rangeF_*$ is connected,  let $\beta$ and $\sigma = F o \beta$ be geodesics on $M$ and $N,$ respectively. Then $F$ is a Clairaut Riemannian map with $\tilde{s} = e^g$ if and only if any one of the following conditions holds: 
	\end{theorem}
	\begin{itemize}
		\item[(i)] $S_D F_*X = -D(g) F_*X$, where $F_*X \in\Gamma (rangeF_*)$ and 
		$D \in\Gamma(rangeF_*)^\perp$ are components of  $\dot{\sigma}(t)$. 
		\item [(ii)] $F$ is an umbilical map, and has $H' = -\nabla^N g$ , where $g$ is a smooth function on $N$ and $H'$ is the mean curvature vector field of $rangeF_*$.
	\end{itemize} 
\begin{definition}
		\cite{S4} Let $F : (M, g_M) \rightarrow (N, g_N, J')$ be a proper Riemannian map from a Riemannian manifold to an almost Hermitian manifold with almost complex structure $J'.$ We say that $F$ is an anti-invariant Riemannian map at $q\in N$ if $J'
		(rangeF_*q) \subset (rangeF_{*q})^\perp.$
	\end{definition}
 If $F$ is an anti-invariant Riemannian map for every $q \in N,$ then $F$ is called an anti-invariant Riemannian map.
	In this case, we denote the orthogonal subbundle to $J'(rangeF_*)$ in $(rangeF_*)^\perp$
	by $\nu$, i.e., $(rangeF_*)^\perp = J'(rangeF_*) \oplus \nu.$\\ For any $D \in \Gamma(rangeF_*)^\perp,$ we have \cite{S50}
	\begin{equation}\label{JV}
		J'D = PD + QD,
	\end{equation} 
	where $PD \in \Gamma (rangeF_*)$ and $QD \in \Gamma (\nu)$. Note that if $\nu = 0$ then $F$ is called a Lagrangian Riemannian map.\\
	Now, we recall the definition of Clairaut anti-invariant Riemannian map for target manifold.  
	\begin{definition}
\cite{S35} An anti-invariant Riemannian map from a Riemannian manifold to a K\"ahler manifold is called Clairaut anti-invariant Riemannian map if it satisfies the definition of Clairaut Riemannian manifold for target manifold.	
	\end{definition}
	\section{Clairaut anti-invariant Riemannian maps from K\"ahler manifolds admitting Ricci soliton to Riemannian manifolds }
	In this section, we study Clairaut anti-invariant Riemannian maps from K\"ahler manifolds admitting Ricci solitons to Riemannian manifolds.
	Throughout this section, we are taking $dim(kerF_*)>1.$ Now, we improve the Proposition (3.9) \cite{S20} to explore the study in our case. 
	\begin{proposition}
		Let $F:(M, g_M, J)\rightarrow (N, g_N)$ be a Clairaut anti-invariant Riemannian map from a K\"ahler manifold to a Riemannian manifold with $\tilde{r} = e^f.$ Then, Ricci tensors of $(M, g_M, J)$ are given by
		\begin{equation}\label{Ric U,V}
			\begin{split}
				Ric(U, V) = Ric^{(range F_*)}(F_*JU, F_*JV)+ rg_M(\nabla_{JU}\nabla f, JV)- div\Big(A(JU, JV)\Big)  
			\end{split}
		\end{equation}
		\begin{equation}\label{Ric U, X}
			\begin{split}
				Ric(U, X)=&-(r+1)g_M(\nabla_{BX}\nabla f, JU)+ div\Big(A(JU, CX)\Big), \\&-rg_M(\nabla_{JU}\nabla f, CX)+Ric^{rangeF_*}(F_*JU, F_*CX)\\&+\sum_{i=1}^{r+s}g_M\Big((\nabla_{X_i}A)(X_i, JU), BX)
			\end{split}	
		\end{equation}
		and 
		\begin{equation}\label{Ric X,Y}
			\begin{split}
				Ric(X, Y)& = Ric^{(KerF_*)}(BX, BY)- \Big(r||\nabla f||^2+div(\nabla f)\Big)g_M(BX, BY)\\&+\sum_{i=1}^{r+s}g_M(A_{X_i}BX, A_{X_i}BY) -rg_M(\nabla_{CX}\nabla f, CY)\\&-rCX(f)CY(f)+\sum_{j=1}^{r}g_M(A_{CX}u_j, A_{CY}u_j)\\&+div\Big(A(CX, CY)\Big)+Ric^{(rangeF_*)}(F_*CX, F_*CY)\\&
				-\sum_{i=1}^{r+s}g_N\Big((\nabla F_*)(X_i, CY), (\nabla F_*)(CX, X_i)\Big)+g_N\Big((\nabla F_*)(CX, CY), \tau^{(kerF_*)^\perp}\Big)\\& 
				-(r+1)g_M(\nabla_{BX}\nabla f, CY) -\sum_{i=1}^{r+s}g_M\Big((\nabla_{X_i}A)(CY, X_i), BX\Big)\\&-(r+1)g_M(\nabla_{BY}\nabla f, CX)-\sum_{i=1}^{r+s}g_M\Big((\nabla_{X_i}A)(CX, X_i), BY\Big) 
			\end{split}
		\end{equation}
		for all $X, Y\in\Gamma(kerF_*)^\perp, U, V\in\Gamma(kerF_*),$ where $\{u_1, u_2,...,u_r\}, \{X_1, X_2,...,X_{r+s}\}$ and $\{\mu_1, \mu_2,...,\mu_s\}$ are orthonormal frames of $kerF_*, J(kerF_*)\oplus \mu$ and $\mu,$ respectively, and Ric, $Ric^{(kerF_*)}$ and $Ric^{rangeF_*}$ are Ricci tensors on $M$, $kerF_*$ and $rangeF_*$, respectively.
	\end{proposition}
	\begin{proof} Let $F$ be an anti-invariant Riemannian map from a K\"ahler manifold $(M, g_M, J)$ to a Riemannian manifold $(N, g_N).$ In \cite{S32}, \c{S}ahin defined $Ric(U, V)$ as
		\begin{equation}\label{Anti Ric(U, V)}
			\begin{split}
				Ric(U, V) = &\sum_{j=1}^{r}\Big\{-g_M\Big((\nabla_{JU}T)(u_j, u_j), JV\Big)-g_M\Big((\nabla_{u_j}A)(JU, JV), u_j\Big)\\&+g_M(T_{u_j}JU, T_{u_j}JV)-g_M(A_{JU} u_j, A_{JV} u_j)\Big\}\\&+Ric^{(rangeF_*)}(F_{*}JU, F_{*}JV)+g_N\Big((\nabla F_*)(JU, JV), \tau^{(kerF_*)^\perp}\Big)\\&-\sum_{i=1}^{r+s}g_N\Big((\nabla F_*)(X_i, JV), (\nabla F_*)(JU, X_i)\Big).		
			\end{split}
		\end{equation}
		Since $F$ is a Clairaut anti-invariant Riemannian map, using Clairaut condition $T_U V = -g(U, V)\nabla f,$ Lemmas 3.7 and  3.8 of \cite{S20}, we get
		$$-\sum_{j=1}^{r}g_M\Big((\nabla_{JU}T)(u_j, u_j), JV\Big) = rg_M(\nabla_{JU}\nabla f, JV),$$
		$$\sum_{j=1}^{r}g_M\Big((\nabla_{u_j}A)(JU, JV), u_j\Big) = div\Big(A(JU, JV)\Big),$$
		$$\sum_{j=1}^{r}g_M(T_{u_j}JU, T_{u_j}JV) = 0,$$
		$$\sum_{j=1}^{r}g_M(A_{JU} u_j, A_{JV} u_j) = 0,$$
		$$g_N\Big(\nabla F_*)(JU, JV), \tau^{(kerF_*)^\perp}\Big) = 0,$$
		$$\sum_{i=1}^{r+s}g_N\Big((\nabla F_*)(X_i, JV), (\nabla F_*)(JU, X_i)\Big) = 0. $$
		Therefore using above results in \eqref{Anti Ric(U, V)}, we get\\
		$$Ric(U, V) = Ric^{(range F_*)}(F_*JU, F_*JV)+ rg_M(\nabla_{JU}\nabla f, JV)- div\Big(A(JU, JV)\Big)\textcolor{red}{,}$$
		Next, in \cite{S32}, $Ric(U, X)$ is defined as
		\begin{equation}\label{Anti Ric(U, X)}
			\begin{split}
				Ric(U, X) = &\sum_{j=1}^{r}\Big\{g_M\Big((\nabla_{BX}T)(u_j, u_j), JU\Big)-g_M\Big((\nabla_{u_j}T)(BX, u_j), JU\Big)\\&+g_M\Big((\nabla_{JU}T)(u_j, u_j), CX\Big) +g_M\Big((\nabla_{u_j}A)(JU, CX), u_j\Big)\\&-g_M(T_{u_j}JU, T_{u_j}CX)+g_M(A_{JU}u_j, A_{CX}u_j)\Big\}+Ric^{(rangeF_*)}(F_*JU, F_*CX)\\&+g_N\Big((\nabla F_*)(JU, CX), \tau^{(kerF_*)^\perp}\Big)+\sum_{i=1}^{r+s}\Big\{g_M\Big((\nabla_{X_i}A)(X_i, JU), BX\Big)\\& +2g_M(A_{X_i}JU, T_{BX}X_i) - g_N\Big((\nabla F_*)(X_i, CX), (\nabla F_*)(JU, X_i)\Big\}.
			\end{split}
		\end{equation}
		Since $F$ is a Clairaut anti-invariant Riemanian map, using Clairaut condition $T_U V = -g(U, V)\nabla f,$ Lemmas 3.7 and 3.8 of \cite{S20}, we get
		\begin{equation*}
			\begin{split}
				&\sum_{j=1}^{r}\Big\{g_M\Big((\nabla_{BX}T)(u_j, u_j), JU\Big)-g_M\Big((\nabla_{u_j}T)(BX, u_j), JU\Big)+2g_M(A_{X_i}JU, T_{BX}X_i)\Big\}\\& = -(r+1)g_M(\nabla_{BX}\nabla f, JU),
			\end{split}
		\end{equation*}
		$$\sum_{j=1}^{r}g_M\Big((\nabla_{u_j}A)(JU, CX), u_j\Big) = div\Big(A(JU, CX)\Big),$$
		$$\sum_{j=1}^{r}g_M(T_{u_j}JU, T_{u_j}CX) = 0,$$
		$$\sum_{j=1}^{r}g_M(A_{JU}u_j, A_{CX}u_j) = 0,$$
		$$g_N\Big((\nabla F_*)(JU, X_i), (\nabla F_*)(X_i, CX) = 0,$$ 
		$$g_N\Big((\nabla F_*)(JU, CX), \tau^{(kerF_*)^\perp}\Big) = 0,$$
		$$g_M\Big((\nabla_{JU}T)(u_j, u_j), CX\Big) = -rg_M(\nabla_{JU}\nabla f, CX).$$
		Using above results in \eqref{Anti Ric(U, X)}, we get
		\begin{equation*}
			\begin{split}
				Ric(U, X)=&-(r+1)g_M(\nabla_{BX}\nabla f, JU)+ div\Big(A(JU, CX)\Big) \\&-rg_M(\nabla_{JU}\nabla f, CX)+Ric^{rangeF_*}(F_*JU, F_*CX)\\&+\sum_{i=1}^{r+s}g_M\Big((\nabla_{X_i}A)(X_i, JU), BX).
			\end{split}	
		\end{equation*}
		Next, since $M$ is a K\"ahler manifold, we have
		\begin{equation}\label{Ric(X, Y) = Ric(JX, JY)}
			Ric(X', Y') = Ric(JX', JY') ~~\forall~~ X', Y'\in \Gamma(TM).	
		\end{equation}
		Using \eqref{JX}, we have\\
		\begin{equation}\label{Ric(JX, JY)}
			Ric(JX, JY) = Ric(BX, BY) + Ric(CX, CY) + Ric(BX, CY) + Ric(CX, BY)
			~~ \forall X, Y\in\Gamma(kerF_*)^\perp.
		\end{equation}
		Next, since
		\begin{equation*}
			Ric(BX, BY) = \sum_{j=1}^{r}g_M(R(u_j, BX)BY, u_j) + \sum_{i=1}^{r+s}g_M(R(X_i, BX)BY, X_i),
		\end{equation*} 
		using Lemma 52 of \cite{S23}, we get
		\begin{equation*}
			\begin{split}
				Ric(BX, BY)& = Ric^{kerF_*}(BX, BY) -rg_M(H, T_{BX}BY) + g_M(A_{X_i}BX, A_{X_i}BY)\\& + \sum_{i=1}^{r+s}g_M\Big((\nabla_{X_i}T)_{BX}BY, X_i \Big). 		
			\end{split}
		\end{equation*}
		Since $F$ is a Clairaut anti-invariant Riemannian map, applying the  Clairaut condition $T_U V = -g_M(U, V)\nabla f$, and metric compatibility condition of the Levi-Civita connection $\nabla$, we get\\
		\begin{equation}\label{Ric(BX, BY)}
			\begin{split}
				Ric(BX, BY)& = Ric^{kerF_*}(BX, BY) - r||\nabla f||^2g_M(BX, BY) + \sum_{i=1}^{r+s}g_M(A_{X_i}BX, A_{X_i}BY)\\& - g_M(BX, BY)div(\nabla f).			
			\end{split}		
		\end{equation} 
		Similar computations provide the following three equations:
		\begin{equation}\label{Ric(CX, CY)}
			\begin{split}
				Ric(CX, CY) =& -rg_M(\nabla_{CX}\nabla f, CY)+Ric^{(rangeF_*)}(F_*CX, F_*CY)\\&-rCX(f)CY(f)+\sum_{j=1}^{r}g_M(A_{CX}u_j, A_{CY}u_j)+div\Big(A(CX, CY)\Big)\\&-\sum_{i=1}^{r+s}g_N\Big((\nabla F_*)(X_i, CY), (\nabla F_*)(CX, X_i)\Big)\\&+g_N\Big((\nabla F_*)(CX, CY), \tau^{(kerF_*)^\perp}\Big),		 
			\end{split}
		\end{equation}
		\begin{equation}\label{Ric(BX, CY)}
			\begin{split}
				Ric(BX, CY) = -(r+1)g_M(\nabla_{BX}\nabla f, CY)-\sum_{i=1}^{r+s} g_M\Big((\nabla_{X_i}A)_{CY}{X_i}, BX),
			\end{split}
		\end{equation}
		\begin{equation}\label{Ric(CX,BY)}
			\begin{split}
				Ric(CX, BY)& = -(r+1)g_M(\nabla_{BY}\nabla f, CX)-\sum_{i=1}^{r+s} g_M\Big((\nabla_{X_i}A)_{CX}{X_i}, BY).
			\end{split}
		\end{equation}
		Adding \eqref{Ric(BX, BY)}, \eqref{Ric(CX, CY)}, \eqref{Ric(BX, CY)}, \eqref{Ric(CX,BY)}, we obtain \eqref{Ric(JX, JY)} and with the help of \eqref{Ric(X, Y) = Ric(JX, JY)}, we get \eqref{Ric X,Y}.
	\end{proof}
	\begin{corollary}
		Let $F:(M,g_M, J)\rightarrow(N, g_N)$ be a Clairaut Lagrangian Riemannian map from a K\"ahler manifold to a Riemannian manifold  with $\tilde{r} = e^f.$ Then, Ricci tensors on $(M, g_M, J)$ are given by
		\begin{equation}
			\begin{split}\label{LRic(U, V)}
				Ric(U, V) =  Ric^{(range F_*)}(F_*JU, F_*JV),		
			\end{split}
		\end{equation}
		\begin{equation}\label{LRic(U, X)}
			\begin{split}
				Ric(U, X)= 0,
			\end{split}	
		\end{equation}
		and
		\begin{equation}\label{LRic(X, Y)}
			\begin{split}
				Ric(X, Y) = Ric^{(KerF_*)}(BX, BY).
			\end{split}
		\end{equation}
	\end{corollary} 
	\begin{proof}
		In case of Lagrangian Riemannian map, $\mu = 0$ or $CX = 0= CY$ and from Lemma 56 of \cite{S32}, we have $A_X = 0~~\forall X\in\Gamma(kerF_*)^\perp.$ From  Theorem 3.5 of \cite{S20}, we get $f$ is constant on $(kerF_*)^\perp.$\\
		Therefore \eqref{Ric U,V}, \eqref{Ric U, X} and \eqref{Ric X,Y} reduce to \eqref{LRic(U, V)}, \eqref{LRic(U, X)} and \eqref{LRic(X, Y)}, respectively. This completes the proof.
	\end{proof}

	\begin{corollary}\label{cor.}
		Let $F:(M, g_M, J)\rightarrow(N, g_N)$ be a totally geodesic Clairaut anti-invariant Riemannian map from a K\"ahler manifold to a Riemannian manifold with with $\tilde{r} = e^f$ such that $(kerF_*)^\perp$ is totally geodesic. Then
		\begin{equation}
			\begin{split}\label{Cor. Ric(X, Y)}
				Ric(X, Y)& = Ric^{(kerF_*)}(BX, BY) -\Big(r||\nabla f||^2+div(\nabla f)\Big)g_M(BX, BY)\\&-rg_M(\nabla_{CX}\nabla f, CY) + Ric^{(rangeF_*)}(F_*CX, F_*CY)-rCX(f)CY(f),		
			\end{split}
		\end{equation}
		where $X, Y\in\Gamma{(kerF_*)^\perp}.$
	\end{corollary}
	\begin{proof}
		Since $F$ is a totally geodesic Clairaut anti-invariant Riemannian map, the second fundamental form of the map $\nabla F_*$ vanishes. Since $(kerF_*)^\perp$ is totally geodesic, we have $A_X = 0 ~~\forall X\in\Gamma(kerF_*)^\perp$. Therefore the \eqref{Ric X,Y}, reduces to \eqref{Cor. Ric(X, Y)}.
	\end{proof}
	\begin{theorem}\label{F(JkerF) RS}
		Let $F:(M, g_M, J)\rightarrow (N, g_N)$ be a Clairaut anti-invariant Riemannian map from a K\"ahler manifold admitting Ricci soliton with potential vector field $\eta\in\Gamma(TM)$ and scalar $\lambda$ to a Riemannian manifold  such that  $(KerF_*)^\perp$ is totally geodesic. Then $F_*(JkerF_*)$ is an $\alpha$- Ricci soliton.
	\end{theorem}
	\begin{proof}
		Since $(M, g_M, J)$ is a K\"ahler manifold admitting a Ricci soliton with potential vector field $\eta\in\Gamma(TM)$ and scalar $\lambda,$ we have\\
		$$\frac{1}{2}(L_\eta g_M)(U, V) +Ric(U, V) + \lambda g_M(U, V) = 0~~\forall~~  U, V\in\Gamma(KerF_*),$$
		which can be rewritten as
		$$\frac{1}{2}\Big(g_M(\nabla_U \eta, V)+ g_M(\nabla_V \eta, U)\Big) +Ric(U, V) + \lambda g_M(U, V) = 0.$$
		Making use of \eqref{Ric U,V} in the above equation, we get\\
		\begin{equation}\label{JkerF}
			\begin{split}
				&\frac{1}{2}\Big(g_M(\nabla_U \eta, V) + g_M(\nabla_V \eta, U)\Big) +rg_M(\nabla_{JU}\nabla f, JV)\\&- div\Big(A(JU, JV)\Big)+ Ric^{(range F_*)}(F_*JU, F_*JV)\\& + \lambda g_M(U, V) = 0.
			\end{split}
		\end{equation} 
		Since 
		$(kerF_*)^\perp$ is totally geodesic, $A_X = 0 ~~ \forall X\in\Gamma(kerF_*)^\perp.$ Thus \eqref{JkerF} reduces to\\
		\begin{equation}\label{1st th.}
			\begin{split}
				&\frac{1}{2}\Big(g_M(\nabla_U \eta, V) + g_M(\nabla_V \eta, U)\Big)+rg_M(\nabla_{JU}\nabla f, JV) + Ric^{(range F_*)}(F_*JU, F_*JV)\\&+ \lambda g_M(U, V) = 0.
			\end{split}
		\end{equation}
		Using \eqref{COHM} and \eqref{COKM} in \eqref{1st th.}, we get\\
		\begin{equation}\label{eqn before eta}
			\begin{split}
				&\frac{1}{2}\Big(g_M(\nabla_U J\eta, JV) + g_M(\nabla_V J\eta, JU)\Big)+rg_M(\nabla_{JU}\nabla f, JV)\\& + Ric^{(range F_*)}(F_*JU, F_*JV)+ \lambda g_M(JU, JV) = 0.		
			\end{split}
		\end{equation}
\textbf{Case (i):}   When $\eta \in\Gamma(kerF_*)$ or $\eta \in \Gamma(\mu).$\\
	In this case \eqref{eqn before eta}
		\begin{equation*}
			\begin{split}
				&\frac{1}{2}\Big(g(A_{J\eta}U, JV)+g(JU, A_{J\eta}V)\Big)+rg_M(\nabla_{JU}\nabla f, JV)\\& + Ric^{(range F_*)}(F_*JU, F_*JV)+ \lambda g_M(JU, JV) = 0.		
			\end{split}
		\end{equation*}
		Since $(kerF_*)^\perp$ is totally geodesic, $A_{J\eta} =0.$ Therefore, above equation gives\\
		\begin{equation*}
			\begin{split}
				&\frac{r}{2}g_N(\nabla^N_{F_*JU}F_*(\nabla f), F_*JV)+\frac{r}{2}g_N(\nabla^N_{F_*JV}F_*(\nabla f), F_*JU)\\& + Ric^{(range F_*)}(F_*JU, F_*JV)+ \lambda g_N(F_*JU, F_*JV) = 0.
			\end{split}
		\end{equation*}
		This gives
		\begin{equation*}
			\begin{split}
				&\frac{1}{2}\Big(g_N(\nabla^N_{F_*JU}F_*(\nabla f), F_*JV)+g_N(\nabla^N_{F_*JV}F_*(\nabla f), F_*JU)\Big)\\& +\alpha Ric^{(rangeF_*)}(F_*JU, F_*JV) +\beta g_N(F_*JU, F_*JV) = 0.
			\end{split}
		\end{equation*}
		Therefore $F_*(JkerF_*)$ is an $\alpha$- Ricci soliton with potential vector field $F_*(\nabla f),$ where $\alpha = \frac{1}{r}$ and $\beta = \frac{\lambda}{r}.$\\
	\textbf{Case (ii):} When $\eta\in\Gamma(KerF_*)^\perp$ but $\notin \Gamma(\mu).$\\
	In this case from \eqref{eqn before eta}, we get\\
		\begin{equation*}
			\begin{split}
				&\frac{1}{2}\Big(g(T_U J\eta, JV)+ g(JU, T_V J\eta)\Big)+\frac{r}{2}g_M(\nabla_{JU}\nabla f, JV)+\frac{r}{2}g_M(\nabla_{JV}\nabla f, JU)\\& + Ric^{(range F_*)}(F_*JU, F_*JV)+ \lambda g_M(JU, JV) = 0.		
			\end{split}
		\end{equation*}
		Using Theorem \ref{CCRmp} and Theorem 3.5 of \cite{S20} , we get\\
		\begin{equation*}
			\begin{split}
				&\frac{1}{2}\Big(g_N(\nabla^N_{F_*JU}F_*(\nabla f), F_*JV)+g_N(\nabla^N_{F_*JV}F_*(\nabla f), F_*JU)\Big)\\& +\alpha Ric^{(rangeF_*)}(F_*JU, F_*JV) +\beta g_N(F_*JU, F_*JV) = 0.
			\end{split}
		\end{equation*}
		This also shows that $F_*(JkerF_*)$ is an $\alpha$- Ricci soliton with potential vector field $F_*(\nabla f)$, where $\alpha = \frac{1}{r}$ and $\beta = \frac{\lambda}{r}.$\\
		Therefore, in either case $F_*(JkerF_*)$ is an $\alpha$-Ricci soliton.
	\end{proof}
	\begin{corollary}
		Let $F$ be a Clairaut Lagrangian Riemannian map from a K\"ahler manifold $(M, g_M, J)$ admitting Ricci soliton with potential vector field $\eta\in(TM)$ and scalar $\lambda$ to a Riemannian manifold $(N, g_N)$ with $\tilde{r}=e^f$. Then $rangeF_*$ is an Einstein.
	\end{corollary}
	\begin{proof}
		Using Theorem \ref{F(JkerF) RS}, Lemma 56 of \cite{S32} and Theorem 3.5 of (\cite{S20}), we get the required result.
	\end{proof}
	\begin{theorem}
		Let $F$ be a Clairaut anti-invariant Riemannian map from a K\"ahler manifold $(M, g_M, J)$ admitting Ricci soliton with potential vector field $\eta\in\Gamma (TM)$ and scalar $\lambda$ to an Einstein manifold $(N, g_N)$ such that $(kerF_*)^\perp$ is totally geodesic. Then $\nabla f$ is a conformal vector field on $J(kerF_*).$
	\end{theorem}
	\begin{proof}
		Let $F$ be a Clairaut anti-invariant Riemannian map from a K\"ahler manifold $(M, g_M, J)$ admitting Ricci soliton with potential vector field $\eta\in\Gamma (TM)$ to an Einstein manifold $(N, g_N),$ then we have $L_\eta = 0$(as we have seen in theorem \ref{F(JkerF) RS}). Then from \eqref{1st th.}, we get
		\begin{equation}\label{conformalVF}
			\begin{split}
				&\frac{r}{2}\Big(g_M(\nabla_{JU}\nabla f, JV)+g_M(\nabla_{JV}\nabla f, JU)\Big)\\& + Ric^{(range F_*)}(F_*JU, F_*JV)+ \lambda g_M(JU, JV) = 0.		
			\end{split}	
		\end{equation}
		Since $rangeF_*$ is an Einstein, $Ric^{(range F_*)}(F_*JU, F_*JV) = \lambda g_N(F_*JU, F_*JV).$\\
		Then from \eqref{N1} and \eqref{conformalVF}, we obtain\\
		\begin{equation*}
			\frac{r}{2}\Big(g_M(\nabla_{JU}\nabla f, JV)+g_M(\nabla_{JV}\nabla f, JU)\Big)	+2\lambda g_M(JU, JV) = 0.
		\end{equation*}
		\begin{equation*}
			\frac{1}{2}(L_{\nabla f}g_M)(JU, JV)+\mu' g_M(JU, JV) = 0,
		\end{equation*}
		where $\mu' = \frac{2 \lambda}{r}$ is constant on $M.$
		This shows that $\nabla f$ is a conformal vector field on $J(kerF_*).$
	\end{proof}
	\begin{theorem}
		Let $F:(M, g_M, J)\rightarrow (N, g_N)$ be a Clairaut Lagrangian Riemannian map from a K\"ahler manifold admitting Ricci soliton with potential vector fields $\eta\in\Gamma(TM)$ and scalar $\lambda$ to a Riemannian manifold. Then $s^{(rangeF_*)} = -\lambda r.$
	\end{theorem}
	\begin{proof}
		Since $(M, g_M, J)$ is a K\"ahler manifold admitting a Ricci soliton with potential vector field $\eta\in\Gamma(TM),$ we have
		\begin{equation*}
			\frac{1}{2}(L_\eta g_M)(U, V) + Ric(U, V) + \lambda g_M(U, V) = 0
			~~\forall~~ U, V\in\Gamma(KerF_*).\\
		\end{equation*}	
	Using Lemma 56 of \cite{S23} and Theorem \ref{F(JkerF) RS}, we get $(L_\eta g_M)(U, V) = 0.$ Then above equation reduces to
		$$Ric(U, V)+ \lambda g_M(U, V) = 0.$$
		Using \eqref{LRic(U, V)} in the above equation, we get
		\begin{equation*}
			\begin{split}
			 Ric^{(rangeF_*)}(F_*JU, F_*JV)+\lambda g_M(U, V)= 0.	
			\end{split}	
		\end{equation*}
		Taking trace of above equation, we get\\
		$$s^{(range F_*)} +\lambda r  = 0,$$
		This completes the proof.
	\end{proof}
	\begin{theorem}
		Let $F:(M,g_M, J)\rightarrow(N, g_N)$ be a Clairaut Lagrangian Riemannian map from a K\"ahler manifold admitting Ricci soliton with vertical potential vector field $\eta\in\Gamma(TM)$ and scalar $\lambda$ to a Riemannian manifold. Then $kerF_*$ is an Einstein.
	\end{theorem}
	\begin{proof}
		Since $(M, g_M, J)$ is a K\"ahler manifold admitting a Ricci soliton with vertical potential vector field $\eta\in\Gamma(TM)$, we have\\
		$$\dfrac{1}{2}(L_{\eta}g_M)(X, Y) + Ric(X, Y) + \lambda g_M(X, Y) = 0.$$
		Using \eqref{LRic(X, Y)} in above equation, we get
		\begin{equation*}
			\begin{split}
				\dfrac{1}{2}(L_{\eta}g_M)(X, Y)+Ric^{(KerF_*)}(BX, BY) + \lambda g_M(X, Y) = 0,
			\end{split}		
		\end{equation*}
		which can be rewritten as
		\begin{equation*}
			\begin{split}
				\dfrac{1}{2}\Big(g_M(\nabla_X {\eta}, Y)+g_M(\nabla_Y {\eta}, X)\Big)+Ric^{(KerF_*)}(BX, BY)+ \lambda g_M(JX, JY) = 0.		
			\end{split}
		\end{equation*}
		Since $\eta$ is vertical potential vector field, above equation gives\\
		\begin{equation}\label{Ric(kerf{BX, BY})}
			Ric^{(KerF_*)}(BX, BY)+ \lambda g_M(BX, BY) = 0.	
		\end{equation}
		 This implies that $kerF_*$ is an Einstein.\\
	\end{proof}
	\begin{corollary}
		Let $F:(M,g_M, J)\rightarrow(N, g_N)$ be a Clairaut Lagrangian Riemannian map from a K\"ahler manifold  admitting Ricci soliton with vertical potential vector field $\eta\in\Gamma(TM)$ and scalar $\lambda$ to a Riemannian manifold. Then scalar curvature of $kerF_*$ is $s^{kerF_*} = -r\lambda.$
	\end{corollary}
	\begin{proof}
		Taking trace of \eqref{Ric(kerf{BX, BY})}, we get the result .	
	\end{proof}
Next, we give a theorem which shows the relation between Ricci curvature and Riemannian metric.
	\begin{theorem}
		Let $F:(M, g_M, J)\rightarrow (N, g_N)$ be a totally geodesic Clairaut anti-invariant Riemannian map from a K\"ahler manifold admitting Ricci soliton with vertical potential vector field $\eta\in\Gamma(TM)$ and scalar $\lambda$ to a Riemannian manifold  such that $(kerF_*)^\perp$ is totally geodesic and $f$ is constant on $\Gamma(\mu).$ Then $$Ric^{(rangeF_*)}(F_*CX, F_*CY) +\lambda g_N(F_*CX, F_*CY) = 0$$ satisfies provided $kerF_*$ is an Einstein.
	\end{theorem}
	\begin{proof}
		Since $(M, g_M, J)$ is a K\"ahler manifold admitting a Ricci soliton with vertical potential vector field $\eta\in\Gamma(TM)$, we have
		\begin{equation*}
			\dfrac{1}{2}(L_{\eta}g_M)(X, Y) + Ric(X, Y) + \lambda g_M(X, Y) = 0
		\end{equation*}
		which can be rewritten as 
		\begin{equation*}
			\dfrac{1}{2}\Big(g_M(\nabla_X\eta, Y)+ g_M(\nabla_Y\eta, X)\Big) + Ric(X, Y) + \lambda g_M(X, Y) = 0.
		\end{equation*}
		Using Corollary \ref{cor.}  in above equation, we get
		\begin{equation*}
			\begin{split}
				&\dfrac{1}{2}\Big(g_M(\nabla_X\eta, Y)+ g_M(\nabla_Y\eta, X)\Big) + Ric^{(kerF_*)}(BX, BY)-rg_M(\nabla_{CX}\nabla f, CY)\\&-(r||\nabla f||^2+div(\nabla f)-\lambda)g_M(BX, BY)+ Ric^{(rangeF_*)}(F_*CX, F_*CY)\\& -rCX(f)CY(f) +\lambda g_N(F_*CX, F_*CY) = 0.				
			\end{split}
		\end{equation*}
		Since $(kerF_*)^\perp$ is totally geodesic and $f$ is constant on $\Gamma(\mu),$  from above equation, we get
		\begin{equation*}
			\begin{split}
				&Ric^{(kerF_*)}(BX, BY)+\lambda g_M(BX, BY)+Ric^{(rangeF_*)}(F_*CX, F_*CY) \\&+\lambda g_N(F_*CX, F_*CY) = 0.			
			\end{split}
		\end{equation*}
		Since $kerF_*$  is an Einstein, then above equation reduces to
		\begin{equation*}
			\begin{split}
			Ric^{(rangeF_*)}(F_*CX, F_*CY) +\lambda g_N(F_*CX, F_*CY) = 0.		
			\end{split}
		\end{equation*}
	This completes the proof.
			\end{proof}
Next, we give a theorem which shows the relation between Ricci curvature and Lie derivative.
	\begin{theorem}
		Let $F:(M, g_M, J)\rightarrow (N, g_N)$ be a Clairaut anti-invariant Riemannian map from a K\"ahler manifold admitting a Ricci soliton with horizontal potential vector field $\eta \in\Gamma(TM)$ and scalar $\lambda$ to a Riemannian manifold such that $(kerF_*)^\perp$ is totally geodesic. Then
		$$Ric^{rangeF_*}(F_*JU, F_*CX) = \frac{1}{2\alpha'}(L_{F_*(\nabla f)}g_N)(F_*JU, F_*CX)$$
	\end{theorem}
	where $\alpha' = \frac{1}{r}.$
	\begin{proof}
		Since $(M, g_M, J)$ is a K\"ahler manifold admitting a Ricci soliton with horizontal potential vector field $\eta\in\Gamma(TM)$, we have
		\begin{equation*}
			\frac{1}{2}(L_{\eta}g_M)(U, X) + Ric(U, X)+ \lambda g_M(U, X) = 0.
		\end{equation*}
		Making use of \eqref{Ric U, X}, we get
		\begin{equation*}
			\begin{split}
				&\frac{1}{2}\Big(g_M(\nabla_U \eta, X)+g_M(\nabla_X \eta,  U)\Big)  -(r+1)g_M(\nabla_{BX}\nabla f, JU)\\&-rg_M(\nabla_{JU}\nabla f, CX)+ div\Big(A(JU, CX)\Big)+Ric^{rangeF_*}(F_*JU, F_*CX)\\&+\sum_{i=1}^{r+s}	g_M\Big((\nabla_{X_i}A)(X_i, JU), BX)\Big) = 0. 
			\end{split}	
		\end{equation*}
		Since $\eta$ is horizontal and $(kerF_*)^\perp$ is totally geodesic, we get
		\begin{equation*}
			\begin{split}
				-rg_M(\nabla_{JU}\nabla f, CX) + Ric^{rangeF_*}(F_*JU, F_*CX) = 0
			\end{split}
		\end{equation*}
		which can be rewritten as
		\begin{equation*}
			\begin{split}
				&\frac{1}{2}\Big(g_N(\nabla_{F_*JU} F_*(\nabla f), F_*CX) +g_N(\nabla_{ F_*{CX}}F_*(\nabla f), F_*JU)\Big)\\& -\frac{1}{r}Ric^{rangeF_*}(F_*JU, F_*CX) = 0.
			\end{split}
		\end{equation*}
		This gives
		\begin{equation*}
			\frac{1}{2}(L_{F_*(\nabla f)}g_N)(F_*JU, F_*CX)- \alpha' Ric^{rangeF_*}(F_*JU, F_*CX)= 0.
		\end{equation*}
		This proves the required result.
	\end{proof}
	\begin{example}\label{EX.}
		Let $(M, g_M, J)$ be a 6-dimensional K\"ahler manifold, where $M= \{(x_1, x_2, x_3, x_4, x_5, x_6)\in \mathbb{R}^6: x_i\neq 0, i=1, \ldots,6\},$ with Riemannian metric
		$$g_M = e^{-2x_4}dx_1^2+e^{-2x_4}dx_2^2+e^{-2x_4}dx_3^2+dx_4^2+dx_5^2+dx_6^2,$$ and $J(x_1, x_2, x_3, x_4, x_5, x_6) = (-x_6, -x_5, -x_4, -x_3, -x_2, -x_1).$\\
		Let $(N, g_N)$ be a 6-dimensional Riemannian manifold, where $N = \{(y_1, y_2, y_3, y_4, y_5, y_6)\in\R^6: y_j\in\R\forall j = 1, \ldots,6\}$ with Riemannian metric
		$$g_N = dy_1^2+e^{-2y_4}dy_2^2+dy_3^2+dy_4^2+e^{-2y_4}dy_5^2+dy_6^2.$$
		Let $F:(M, g_M, J)\rightarrow (N, g_N)$ be a map defined as
		$$F(x_1, x_2, x_3, x_4, x_5, x_6) = (x_5, x_2, 0, x_4, 0, x_6).$$ 
		Then by direct computation, we get\\
		$$kerF_* = span\Big\{U_1 = e_1, U_2 = e_3\Big\},$$
		and
		$$(kerF_*)^\perp = span\Big\{X_1 = e_2, X_2 = e_4, X_3 = e_5, X_4 = e_6\},$$
		where $$\{e_1 = e^{x_4}\frac{\partial}{\partial x_1}, e_2 = e^{x_4}\frac{\partial}{\partial x_2}, e_3 = e^{x_4}\frac{\partial}{\partial x_3}, e_4 =\frac{\partial}{\partial x_4}, e_5 = \frac{\partial}{\partial x_5}, e_6 = \frac{\partial}{\partial x_6}\}$$ is the orthonormal basis of $T_pM$, for any $p\in M.$\\
		Then, we obtain	\\$F_*(X_1) = e_2',$ $F_*(X_2) =e_4', F_*(X_3) = e_1', F_*(X_4)= e_6',$\\
		where $\{e_1' = \frac{\partial}{\partial y_1}, e_2' = e^{y_4} \frac{\partial}{\partial y_2}, e_3' =\frac{\partial}{\partial y_3}, e_4' =\frac{\partial}{\partial y_4}, e_5'=e^{y_4}\frac{\partial}{\partial y_5}, e_6'= \frac{\partial}{\partial y_6} \}$ represents the standard orthonormal basis of $T_{F(p)}N.$ \\
		Also, it can be easily verified that  for $i=1,...,4$ $$g_M(X_i, X_i) = g_N(F_*{X_i}, F_*{X_i}),$$ and $J(U_1) = X_1, J(U_2) = X_2, J(X_1) = -U_1, J(X_2) = -U_2, J(X_3) = X_4, J(X_4) = -X_3.$ Thus $F$ is an anti-invariant Riemannian map.\\
		\noindent	Now in order to show that $F$ is a Clairaut anti-invariant Riemannian map, we need to find a smooth function $f$ on $M$ such that $$T_U U = -g(U, U)\nabla f ~~\text{for} ~~U\in\Gamma(kerF_*).$$
		For the Riemannian metric $g_M,$  the Christoffel symbols are obtained as\\
		$\Gamma^1_{1 4} = -1=\Gamma^1_{4 1}, \Gamma^4_{1 1} = e^{-2x_4}, \Gamma^4_{2 2}= e^{-2x_4}, \Gamma^4_{3 3}= e^{-2x_4}, \Gamma^2_{2 4}= -1 = \Gamma^2_{4 2}, \Gamma^3_{3 4} = -1 = \Gamma^3_{4 3}$ and rest of them are zero.\\
		\noindent Additionally, the covariant derivatives are calculated as\\
		$\nabla_{U_1}U_1 = X_2, \nabla_{U_1}X_2 = -U_1,\nabla_{X_1}X_1 = X_2, \nabla_{U_2}U_2 = X_2, \nabla_{X_1}X_2 = -X_1, \nabla_{U_2}X_2 = -U_2,$ rest of them are zero.\\
		Ricci curvatures are obtained as
		$Ric(U_1, U_1) = 3, Ric(U_2, U_2) = 1, Ric(U_1, U_2) = -1, Ric(X_1, X_1) = -2, Ric(X_2, X_2) = -2$
		remaining Ricci tensors are zero.\\
		As we know, $$TM = ker F_*\oplus(kerF_*)^\perp,$$ and for i = $1,\ldots6$, we have $$Z_i = a_iU_1+b_iU_2+c_iX_1+d_iX_2+l_iX_3+h_iX_4. $$
		Also, we have $$(L_{Z_1}g_M)(Z_2, Z_3) = g_M(\nabla_{Z_2}Z_1, Z_3) +g_M(\nabla_{Z_3}Z_1, Z_2).$$
		By, utilizing the values of covariant derivatives, we obtain \\
		\begin{equation}
			\begin{split}\label{L for 3.1}
				(L_{Z_1}g_M)(Z_2, Z_3) =& \frac{1}{2}(a_1a_2d_3+a_1a_3d_2+c_3c_1d_2+b_3b_1d_2+c_2c_1d_3+b_2b_1d_3)\\&-a_2a_3d_1-c_3c_2d_1-b_2b_3d_1,
			\end{split}
		\end{equation}
		\begin{equation}\label{R for 3.1}
			Ric(Z_2, Z_3) = 3a_2a_3+b_2b_3-a_2b_3-2c_2c_3-2d_2d_3,
		\end{equation}
		and
		\begin{equation}\label{g for 3.1}
			g_M(Z_2, Z_3) = a_2a_3+b_2b_3+c_2c_3+d_2d_3+l_2l_3+h_2h_3.
		\end{equation}
		Let $(M, g_M, J)$ admits a Ricci soliton, then
		$$\frac{1}{2}(L_{Z_1}g_M)(Z_2, Z_3) + Ric(Z_2, Z_3) +\lambda g_M(Z_2, Z_3) = 0.$$
		Using \eqref{L for 3.1}, \eqref{R for 3.1} and \eqref{g for 3.1} in the above equation, we get\\
		\begin{equation*}
			\begin{split}
				\lambda = \frac{A_1-A_2}{2(a_2a_3+b_2b_3+c_2c_3+d_2d_3+l_2l_3+h_2h_3)},~(a_2a_3+b_2b_3+c_2c_3+d_2d_3+l_2l_3+h_2h_3) \neq 0,
			\end{split}
		\end{equation*} 
		where $$A_1 = 2(3a_2a_3+b_2b_3-a_2b_3-2c_2c_3-2d_2d_3),$$
		$$A_2 = (a_1a_2d_3+a_1a_3d_2+c_3c_1d_2+b_3b_1d_2+c_2c_1d_3+b_2b_1d_3)-2a_2a_3d_1-2c_3c_2d_1-2b_2b_3d_1.$$
		For some particular values of $a_i's, b_i's, c_i's, d_i's, l_i's, \text{and}~ h_i's,$ the behaviour of Ricci soliton will be shrinking, steady, and expanding accorrding to $\lambda<0, \lambda=0,~\text{and} ~\lambda>0$ respectively.
	\end{example}

				\section{ Clairaut Anti-invariant Riemannian maps to K\"ahler manifolds admitting Ricci Solitons}
				In this section, we study Clairaut anti-invariant Riemannian maps from Riemannian manifolds to K\"ahler manifolds admitting Ricci solitons.
				Throughout this section, $(rangeF_*)^\perp$ is assumed to be totally geodesic and $dim(rangeF_*)>1.$
				\begin{theorem}
					Let $F : (M, g_M) \rightarrow (N, g_N, J')$ be a Clairaut anti-invariant
					Riemannian map from a Riemannian manifold to a K\"ahler manifold with $\tilde{s} = e^g$ . Then Ricci tensors on $(N, g_N, J')$ are
					\begin{equation}\label{Ric F_*x F_*Y}
						Ric(F_*X, F_*Y ) = Ric^{(rangeF_*)^\perp}(J'F_*X, J'F_*Y ) + (m-r)g_N(\nabla^Ng, \nabla^{F\perp}_{J'F_*X} J'F_*Y)
					\end{equation}
					\begin{equation}\label{Ric(F_*X E)}
						\begin{split}
							Ric(F_*X, E) =& Ric^{(rangeF_*)^\perp}(J'F_*X, QE) + \sum_{j=r+1}^{m}g_N\Big((\tilde{\nabla}_{^{*}F_{*}PE}S)_{J'F_*X} F_*X_j,  F_*X_j\Big)\\&- \sum_{j=r+1}^{m}g_N\Big((\tilde{\nabla}_{{X_j}}S)_{J'F_*X} PE,  F_*X_j\Big) +(m-r)g_N(\nabla^N g, \nabla^{F\perp}_{QE} J'F_*X)\\& - \sum_{k=1}^{r_1}g_N(R^{F\perp}(PE, e_{k})J'F_*X, e_{k}),
						\end{split}
					\end{equation}
					\begin{equation}\label{RIC(D, E)}
						\begin{split}
							Ric(D, E) =& Ric^{rangeF_*}(PD, PE)-\sum_{k=1}^{r_1}g_N(PD, PE)\Big(||\nabla^N g||^2+H^g(e_k, e_k)\Big)\\&+\sum_{j=r+1}^{m}g_N\Big((\tilde{\nabla}_{^{*}_{F_*PD}}S)_{QE}F_*{X_j}, F_*{X_{j}}\Big)-\sum_{j=r+1}^{m}g_N\Big((\tilde{\nabla}_{X_j}S)_{QE}PD, F_*{X_j}\Big)\\&-\sum_{k=1}^{r_1}g_N(R^{F\perp}(PD, e_k)QE, e_k)+\sum_{j=r+1}^{m}g_N\Big((\tilde{\nabla}_{^{*}_{F_*PE}}S)_{QD}F_*{X_j}, F_*{X_j}\Big)\\&-\sum_{j=r+1}^{m}g_N\Big((\tilde{\nabla}_{X_j}S)_{QD}PE, F_*{X_j}\Big)-\sum_{k=1}^{r_1}g_N(R^{F\perp}(PE, e_k)QD, e_k)\\&+Ric^{(rangeF_*)^\perp}(QD, QE)-(m-r)\Big(QD(g)QE(g)+g_N(\nabla^N_{QD}\nabla^N g, QE)\Big)
						\end{split}
					\end{equation}
					for $X, Y \in \Gamma(kerF_*)^\perp, F_*X, F_*Y \in \Gamma(rangeF_*)$ and $D, E \in \Gamma(rangeF_*)^\perp$, where
					$\{F_*{X_j}\}_{r+1\leq j\leq m}$ and $\{e_k\}_{1\leq k\leq n_1}$ are orthonormal bases of $rangeF_*$ and $(rangeF_*)^\perp$,
					respectively, where $^{*}{F_*}$ is the adjoint map of $F_*$ and $Ric^{(rangeF_*)},$ $Ric^{(rangeF_*)^\perp}$ are Ricci tensors on $(rangeF_*),$ $(rangeF_*)^\perp,$ respectively.		
				\end{theorem}
				\begin{proof}
					Making use of Proposition 5.1 of \cite{S35}, we get \eqref{Ric F_*x F_*Y} and \eqref{Ric(F_*X E)}. For the proof of \eqref{RIC(D, E)}, we proceed as follows.\\
					Since $N$ is a K\"ahler manifold, we have
					\begin{equation}\label{Ric(DE) = Ric(JD, JE)}
						Ric(D, E) = Ric(JD, JE) ~~\forall D, E\in\Gamma(rangeF_*)^\perp.
					\end{equation}
					Using \eqref{JV}, we get
					\begin{equation}\label{Ric(JDJE)}
						Ric(JD, JE) = Ric(PD, PE)+Ric(QD, PE)+Ric(PD, QE)+ Ric(QD, QE).
					\end{equation}
					Since $$Ric(PD, PE) = \sum_{j=r+1}^{m}g_N(R(F_*{X_j}, PD)PE, F_*{X_j}) + \sum_{k=1}^{n_1}g_N(R(e_k, PD)PE, e_k).$$
					Using Theorem 3.1 of \cite{S51} in above equation, we get
					\begin{equation*}
						\begin{split}
							Ric(PD, PE)& = Ric^{(rangeF_*)}(PD, PE)-\sum_{k=1}^{n_1}g_N(S_{\nabla^{F\perp}_{e_k}e_k}PD, PE)+\sum_{k=1}^{n_1}g_N(\nabla^N_{e_k}S_{e_k}PD, PE)\\&-\sum_{k=1}^{n_1}g_N(S_{e_k}PD, S_{e_k}PE)-\sum_{k=1}^{n_1}g_N(\nabla^N_{e_k}PD, S_{e_k}PE).	
						\end{split}	
					\end{equation*}
					Using Theorem 3.1 of \cite{S35} in above equation, we get\\
					\begin{equation*}
						\begin{split}
							Ric(PD, PE)& = Ric^{(rangeF_*)}(PD, PE)+\sum_{k=1}^{n_1}g_N(PD, PE)g_N(\nabla^N g, \nabla^{F\perp}_{e_k}e_k)\\&-\sum_{k=1}^{n_1}g_N(\nabla^N_{e_k}(e_k(g)PD), PE)+\sum_{k=1}^{n_1}g_N(\nabla^N_{e_k}PD, e_k(g)PE)\\&-\sum_{k=1}^{n_1}(e_k(g))^2g_N(PD, PE).	
						\end{split}
					\end{equation*}
					After simplification, we get
					\begin{equation}\label{Ric(PDPE)}
						\begin{split}
							Ric(PD, PE)& = Ric^{(rangeF_*)}(PD, PE)-\sum_{k=1}^{n_1}g_N(PD, PE)\Big(||\nabla^N g||^2+H^g(e_k, e_k)\Big).
						\end{split}
					\end{equation}
					By similar computation, we get
					\begin{equation}\label{Ric(PDQE)}
						\begin{split}
							Ric(PD, QE) =& \sum_{j=r+1}^{m}g_N\Big((\tilde{\nabla}_{^{*}F_{*}PD}S)_{QE}F_*{X_j}, F_*{X_{j}}\Big)-\sum_{j=r+1}^{m}g_N\Big((\tilde{\nabla}_{X_j}S)_{QE}PD, F_*{X_j}\Big)\\&-\sum_{k=1}^{n_1}g_N(R^{F\perp}(PD, e_k)QE, e_k),		
						\end{split}
					\end{equation}
					\begin{equation}\label{Ric(PEQD)}
						\begin{split}
							Ric(PE, QD) = &\sum_{j=r+1}^{m}g_N\Big((\tilde{\nabla}_{^{*}F_{*}PE}S)_{QD}F_*{X_j}, F_*{X_j}\Big)\\&-\sum_{j=r+1}^{m}g_N\Big((\tilde{\nabla}_{X_j}S)_{QD}PE, F_*{X_j}\Big)-\sum_{k=1}^{n_1}g_N(R^{F\perp}(PE, e_k)QD, e_k)
						\end{split}
					\end{equation}
					and
					\begin{equation}\label{Ric(QDQE)}
						\begin{split}
							Ric(QD, QE) = Ric^{(rangeF_*)^\perp}(QD, QE)-(m-r)\Big(QD(g)QE(g)+g_N(\nabla^N_{QD}\nabla^N g, QE)\Big).
						\end{split}
					\end{equation}
					Adding \eqref{Ric(PDPE)}, \eqref{Ric(PDQE)}, \eqref{Ric(PEQD)} and \eqref{Ric(QDQE)}, we get \eqref{Ric(JDJE)} and with the help of \eqref{Ric(DE) = Ric(JD, JE)}, we get \eqref{RIC(D, E)}.
				\end{proof}
				\begin{corollary}
					Let $F : (M, g_M) \rightarrow (N, g_N, J')$ be a Clairaut Lagrangian
					Riemannian map from a Riemannian manifold to a K\"ahler manifold with $\tilde{s} = e^g$. Then Ricci curvature is given by 
					\begin{equation}\label{LRic F_*x F_*Y}
						Ric(F_*X, F_*Y ) = Ric^{(rangeF_*)^\perp}(J'F_*X, J'F_*Y ) 
					\end{equation}
					\begin{equation}\label{LRic(F_*X E)}
						\begin{split}
							Ric(F_*X, E) =&\sum_{j=r+1}^{m}g_N\Big((\tilde{\nabla}_{^{*}_{F_*PE}}S)_{J'F_*X} F_*X_j,  F_*X_j)\\&- \sum_{j=r+1}^{m}g_N\Big((\tilde{\nabla}_{X_j}S)_{J'F_*X} PE,  F_*X_j)\\& -\sum_{k=1}^{r_1} g_N(R^{F\perp}(PE, e_{k})J'F_*X, e_{k}),
						\end{split}
					\end{equation}
					\begin{equation}\label{LRIC(D, E)}
						\begin{split}
							Ric(D, E) = Ric^{(rangeF_*)}(PD, PE)
						\end{split}
					\end{equation}
					where $D, E\in\Gamma({rangeF_*})^\perp.$
				\end{corollary}
				\begin{proof}
					In case of Lagrangian Riemannian map, $\nu = 0,$ this implies $QD = QE = 0.$ From Theorem 5.2 of \cite{S35}, we have $f$  is constant on $J'(rangeF_*).$ Applying these things in \eqref{Ric F_*x F_*Y}, \eqref{Ric(F_*X E)} and \eqref{RIC(D, E)}, we get the required result.
				\end{proof}
				\begin{theorem}\label{J(rang) Eins}
					Let $F:(M, g_M)\rightarrow(N, g_N, J')$ be a totally geodesic Clairaut anti-invariant Riemannian map from a Riemannian manifold  to a K\"ahler manifold admitting Ricci soliton with vertical potential vector field $\vartheta$ and scalar $\lambda$ such that $(rangeF_*)^\perp$ is an Einstein. Then $\vartheta$ is a conformal vector field on $rangeF_*$.
				\end{theorem}
				\begin{proof}
					Since $(N, g_N, J')$ is a K\"ahler manifold admitting a Ricci soliton with vertical potential vector field $\vartheta\in\Gamma(TN)$,  we have\\
					$$\frac{1}{2}(L_\vartheta g_N)(F_*X, F_*Y) +Ric(F_*X, F_*Y) + \lambda g_N(F_*X, F_*Y) = 0.$$
					Making use of \eqref{Ric F_*x F_*Y} in above equation, we obtain\\
					\begin{equation*}
						\begin{split}
							&\frac{1}{2}(L_\vartheta g_N)(F_*X, F_*Y) +Ric^{(rangeF_*)^\perp}(J'F_*X, J'F_*Y)\\& + (m-r)g_N(\nabla^Ng, \nabla^{F\perp}_{J'F_*X} J'F_*Y) + \lambda g_N(F_*X, F_*Y) = 0,			
						\end{split}
					\end{equation*}
					which reduces to 
					\begin{equation*}
						\begin{split}
							&\frac{1}{2}\Big(g_N(\nabla_{F_*X}\vartheta, F_*Y)+g_N(\nabla_{F_*Y}\vartheta, F_*X)\Big) +Ric^{(rangeF_*)^\perp}(J'F_*X, J'F_*Y) \\& + (m-r)g_N(\nabla^Ng, \nabla^{F\perp}_{J'F_*X} J'F_*Y) + \lambda g_N(F_*X, F_*Y) = 0.		
						\end{split}
					\end{equation*}
					Since $F$ is a totally geodesic Clairaut anti-invariant Riemannian map, second fundamental form of $F$ vanishes. Then by Theorem 3.2 of \cite{S35}, g is constant on $N.$ Therefore, we get
					\begin{equation}\label{Ricci of J(range)}
						\begin{split}
							&\frac{1}{2}\Big(g_N(\nabla^N_{F_*X}\vartheta, F_*Y)+g_N(\nabla^N_{F_*Y}\vartheta, F_*X)\Big) + Ric^{(rangeF_*)^\perp}(J'F_*X, J'F_*Y)\\&+ \lambda g_N(F_*X, F_*Y) = 0.		
						\end{split}
					\end{equation}
					Since $\vartheta$ is vertical, i.e., $\vartheta = F_*Z,$ then from above equation, we get\\
					\begin{equation}\label{thm J(range)}
						\begin{split}
							&\frac{1}{2}\Big(g_N(\nabla^N_{F_*X}F_*Z, F_*Y)+g_N(\nabla^N_{F_*Y}F_*Z, F_*X)\Big)+ Ric^{(rangeF_*)^\perp}(J'F_*X, J'F_*Y)\\&+\lambda g_N(J'F_*X, J'F_*Y)= 0.
						\end{split}		
					\end{equation}
					Since $(rangeF_*)^\perp$ is an Einstein, $$Ric^{(rangeF_*)^\perp}(J'F_*X, J'F_*Y) = \lambda g_N(J'F_*X, J'F_*Y) = \lambda g_N(F_*X, F_*Y),$$ 
					using above equation in \eqref{thm J(range)}, we get\\
					\begin{equation*}
						\begin{split}
							&\frac{1}{2}\Big(g_N(\nabla^N_{F_*X}F_*Z, F_*Y)+g_N(\nabla^N_{F_*Y}F_*Z, F_*X)\Big)+2\lambda g_N(F_*X, F_*Y)= 0,
						\end{split}		
					\end{equation*}
					this gives
					\begin{equation*}
						\frac{1}{2}(L_{F_*Z}g_N)(F_*X, F_*Y)+ 2\lambda g_N(F_*X, F_*Y)= 0.
					\end{equation*}
					This implies that $F_*Z$ is a conformal vector field on $(rangeF_*).$\\
				\end{proof}
				\begin{theorem}\label{Horizontal}
					Let $F:(M, g_M)\rightarrow(N, g_N, J')$ be a Clairaut Lagrangian Riemannian map from a Riemannian manifold  to a K\"ahler manifold admitting Ricci soliton with horizontal potential vector field $\vartheta$ and scalar $\lambda.$ Then $(rangeF_*)^\perp$ is an Einstein.
				\end{theorem}
				\begin{proof}
					Since $(N, g_N, J')$ is a K\"ahler manifold admitting a Ricci soliton with horizontal potential vector field $\vartheta\in\Gamma(TN)$,  we have\\
					$$\frac{1}{2}(L_\vartheta g_N)(F_*X, F_*Y) +Ric(F_*X, F_*Y) + \lambda g_N(F_*X, F_*Y) = 0.$$
					Using \eqref{LRic F_*x F_*Y} in above equation, we get
					\begin{equation*}
						\begin{split}
							&\frac{1}{2}(L_\vartheta g_N)(F_*X, F_*Y) +Ric^{(rangeF_*)^\perp}(J'F_*X, J'F_*Y)\\& + \lambda g_N(F_*X, F_*Y) = 0.			
						\end{split}
					\end{equation*}
				 Since $\vartheta$ is horizontal, i.e., $\vartheta = D,$ with the help of \eqref{S_V}, and Theorem 3.1 of \cite{S35}, we get
							\begin{equation}\label{D(g)}
							D(g)g_N(F_*X, F_*Y)  + Ric^{(range F_*)^\perp}(J'F_*X, J'F_*Y) + \lambda g_N(F_*X, F_*Y) = 0.
						\end{equation}
						Since $N$ is a K\"ahler manifold, we have $$g_N(F_*X, F_*Y) = g_N(J'F_*X, J'F_*Y),$$ then applying above equation in \eqref{D(g)}, we get
						\begin{equation}\label{range perp}
							Ric^{(range F_*)^\perp}(J'F_*X, J'F_*Y) + \rho g_M(J'F_*X, J'F_*Y) = 0,
						\end{equation} where $\rho = D(g)+\lambda.$
					This proves that $(rangeF_*)^\perp$ is an Einstein.
				\end{proof}					
				
				\begin{corollary}
					Let $F:(M, g_M)\rightarrow(N, g_N, J')$ be a Clairaut Lagrangian Riemannian map from a Riemannian manifold to a K\"ahler manifold admitting Ricci soliton with horizontal potential vector field $\vartheta\in\Gamma(TN)$ and scalar $\lambda.$ Then scalar curvature $s^{(rangeF_*)^\perp}$ of $(rangeF_*)^\perp$ is given by
					\begin{equation}
						s^{(rangeF_*)^\perp} = -\rho n_1,
					\end{equation}
					where $\rho= D(g)+\lambda.$
				\end{corollary}
					
				\begin{proof}
					Taking trace of \eqref{range perp}, we get the required result.
				\end{proof}
				\begin{theorem}
					Let $F: (M, g_M)\rightarrow (N, g_N, J')$ be a Clairaut Lagrangian Riemannian map from a Riemannian manifold to a K\"ahler manifold admitting Ricci soliton with vertical potential vector field $\vartheta$ and scalar $\lambda.$ Then $rangeF_*$ is an Einstein.
				\end{theorem}
				\begin{proof}
					Since $(N, g_N, J')$ is a K\"ahler manifold admitting a Ricci soliton with vertical potential vector field $\vartheta\in\Gamma(TN)$,  we have\\ $$\frac{1}{2}(L_{\vartheta}g_N)(D, E) + Ric(D, E) + \lambda g_N(D, E) = 0.$$
					Using \eqref{LRIC(D, E)} in above equation, we get\\
					\begin{equation}
						\begin{split}
							&\frac{1}{2}\Big(g_N(\nabla_D \vartheta, E) + g_N(\nabla_E \vartheta, D)+ Ric^{(rangeF_*)}(PD, PE)\\& + \lambda g_N(PD, PE)= 0.		
						\end{split}
					\end{equation}
					Since $\vartheta$ is vertical, i.e., $\vartheta = F_*X,$ we get
					\begin{equation*}
						\begin{split}
							&\frac{1}{2}\Big(g_N(\nabla_D F_*X, E) + g_N(\nabla_E F_*X, D)+ Ric^{(rangeF_*)}(PD, PE)\\& + \lambda g_N(PD, PE) = 0.		
						\end{split}	
					\end{equation*}
					Using metric compatibility, we obtain
					\begin{equation}\label{kappa}
						Ric^{(rangeF_*)}(PD, PE)+ \lambda g_N(PD, PE) = 0,
					\end{equation}
					This shows that $rangeF_*$ is an Einstein.  
				\end{proof}
				\begin{theorem}
					Let $F:(M, g_M)\rightarrow (N, g_N, J')$ be a Clairaut Lagrangian Riemannian map from a Riemannian manifold  to a K\"ahler manifold admitting Ricci soliton with horizontal potential vector field $\vartheta$ and a scalar $\lambda.$ Then the scalar curvature $s^{(rangeF_*)}$ of $(rangeF_*)$ is given by $s^{(rangeF_*)} = -(m-r)\lambda.$ 
				\end{theorem}
				\begin{proof}
					Since $(N, g_N, J')$ is a K\"ahler manifold admitting a Ricci soliton with horizontal potential vector field $\vartheta\in\Gamma(TN)$,  we have\\
					\begin{equation*}
						\begin{split}
							\frac{1}{2}(L_\vartheta g_N)(D, E) + Ric(D, E) +\lambda g_N(D, E) = 0.	 
						\end{split}	
					\end{equation*}	
					Since $\vartheta$ is horizontal vector field and using \eqref{LRIC(D, E)}, we have
					\begin{equation*}
						\begin{split}
							\frac{1}{2}\Big(g_N(\nabla_D \vartheta, E)+ g_N(D, \nabla_E \vartheta)\Big) + Ric^{(range F_*)}(PD, PE)+\lambda g_N(PD, PE) = 0.
						\end{split}	
					\end{equation*}
					Taking trace of above equation, we get
					$$\sum_{k=1}^{r_1}g_N(\nabla_{e_k} e_k, e_k) + \sum_{j=r+1}^{m}Ric^{(range F_*)}(F_*{X_j}, F_*{X_j})+\sum_{j=r+1}^{m}\lambda g_N(F_*{X_j}, F_*{X_j}) = 0$$
					which gives
					\begin{equation}
						s^{(rangeF_*)} = -(m-r)\lambda. 	
					\end{equation}
					If $\vartheta$ is vertical vector field, then taking trace of \eqref{kappa}, we get the same result.
					This completes the proof.
				\end{proof}

				\begin{example}
					Let $M = \{(x_1, x_2, x_3, x_4, x_5, x_6): x_i\neq 0, i= 1,...,6\}$ be a Riemannian manifold with Riemannian metric $g_M = dx_1^2+ dx_2^2+e^{2x_5}dx_3^2+e^{2x_5}dx_4^2+dx_5^2+e^{2x_5}dx_6^2$ and $N = R^6$ be a K\"ahler manifold with Riemannian metric and complex structure $g_N = dy_1^2+dy_2^2+e^{2y_5}dy_3^2+dy_4^2+dy_5^2+e^{2y_5}dy_6^2,$ $J'(a, b, c, d, l, h) = (-b, a, -d, c, -h, l)$ respectively. We construct a map $F: (M, g_M)\rightarrow (N, g_N, J'),$ defined by\\
					$$F(x_1, x_2, x_3, x_4, x_5, x_6) = \big(0, x_2, 0, 0, x_5, 0\Big)$$
					Then $kerF_* = span\{U_1 = e_1, U_2 = e_3, U_3 = e_4, U_4 = e_6\}$ and $(kerF_*)^\perp = span\{X_1 =e_2, X_2 = e_5\},$
					where $\{e_1 = \frac{\partial}{\partial x_1}, e_2 = \frac{\partial}{\partial x_2}, e_3 = e^{-x_5} \frac{\partial}{\partial x_3}, e_4 = e^{-x_5} \frac{\partial}{\partial x_4}, e_5 = \frac{\partial}{\partial x_5}, e_6 = e^{-x_5} \frac{\partial}{\partial x_6}\}$ and $\{e_1' = \frac{\partial}{\partial y_1}, e_2' = \frac{\partial}{\partial y_2}, e_3' = e^{-y_5}\frac{\partial}{\partial y_3}, e_4' = \frac{\partial}{\partial y_4}, e_5' = \frac{\partial}{\partial y_5}, e_6' =e^{-y_5} \frac{\partial}{\partial y_6}\}$ are the bases of $T_pM$ and $T_{F(p)N}$ respectively.\\
					Here $F_*(X_1) = e_2'$ and $F_*(X_2) = e_5'$ and we can easily compute that $$g_M(X, X) = g_N(F_*X, F_*X) ~~\forall~ X\in(kerF_*)^\perp.$$
					$rangeF_* = span\{ F_*{X_1} = e_2', F_*{X_2} = e_5'\}$ and $(rangeF_*)^\perp = \{e_1', e_3', e_4', e_6'\}.$\\
					Also we can see that $J'F_*{X_1} = -e_1', J'F_*{X_2} = e_6'$ which implies that $F$ is an anti-invariant Riemannian map.\\
					To show $F$ is a Clairaut anti-invariant Riemanian map, we have to find $g$ such that $$(\nabla F_*)(X, X) = -g(X, X)\nabla^N g.$$
					After some computation, we get \\
					$\Gamma^3_{3 5} = 1 = \Gamma^3_{5 3}, \Gamma^6_{5 6} = 1 = \Gamma^6_{6 5}, \Gamma^5_{3 3} = -e^{2y_5}, \Gamma^5_{6 6} = -e^{2y_5}$ and rest of Christoffel symbols are zero.\\
					The covariant derivative is given as\\
					$\nabla_{e_3'} {e_5'} = e_3', \nabla_{e_6'}e_5'= e_6', \nabla_{e_3'} {e_3'} = -e_5', \nabla_{e_6'} e_6' = -e_5'$ and rest of covariant derivatives are zero.\\

					Since, $(\nabla F_*)(X,X) \in \Gamma(rangeF_*)^\perp$ for any
					$X \in \Gamma(kerF_*)^\perp.$ So, here we can write $(\nabla
					F_*)(X,X) = \eta_1e_1'+\eta_2e_3'+\eta_3e_4'+\eta_4e_6',$
					and $g_M(X,X) = g_M(\delta_1X_1+\delta_2X_2, \delta_1X_1+\delta_2X_2) = \delta_1^2+\delta_2^2 = \delta^2$ for some $\eta_1, \eta_2, \delta_1, \delta_2 \in \mathbb{R}.$\\		
					Gradient of any smooth function $g$ is given by
					$$\nabla^N g = \sum_{i,j = 1}^{6}g^{ij}\frac{\partial g}{\partial y_i}\frac{\partial}{\partial y_j}.$$ Therefore 
					$$\nabla^N g = \frac{\partial g}{\partial y_1}\frac{\partial}{\partial y_1} + \frac{\partial g}{\partial y_2}\frac{\partial}{\partial y_2} + e^{-2y_5}\frac{\partial g}{\partial y_3}\frac{\partial}{\partial y_3} + \frac{\partial g}{\partial y_4}\frac{\partial}{\partial y_4} + \frac{\partial g}{\partial y_5}\frac{\partial}{\partial y_5} + e^{-2y_5}\frac{\partial g}{\partial y_6}\frac{\partial}{\partial y_6}.$$
					Thus $$\nabla^N g = \frac{\eta_1}{\delta^2}e_1'+ e^{-y_5}\frac{\eta_2}{\delta^2}e_3'+ \frac{\eta_3}{\delta^2}e_4' + e^{-y_5}\frac{\eta_4}{\delta^2}e_6'$$ for the function $g = \frac{\eta_1}{\delta^2}y_1 + e^{-y_5}\frac{\eta_2}{\delta^2}y_3 + \frac{\eta_3}{\delta^2} y_4 + e^{-y_5}\frac{\eta_4}{\delta^2}y_6.$
					Further, proceeding as in example \ref{EX.}, one can easily show that $N$ admits a  Ricci soliton.
				\end{example}
				\section{Conclusion and future scope}
				We have successfully investigated the geometry of Clairaut anti-invariant Riemannian maps with K\"ahler and Ricci soliton structures. In the future, one can study Clairaut invariant \cite{CIRM}, semi-invariant \cite{kpk_tjm, Polat-Meena} and slant \cite{S52} Riemannian maps with K\"ahler and Ricci soliton structures.
				
				\bigskip
				
				\section*{Acknowledgement}
				
				First author is grateful to the financial support provided by CSIR (Council
				of science and industrial research) Delhi, India. File
				no.[09/1051(12062)/2021-EMR-I]. The second author is thankful to the
				Department of Science and Technology(DST) Government of India for providing
				financial assistance in terms of FIST project(TPN-69301) vide the letter
				with Ref No.:(SR/FST/MS-1/2021/104).\\ 

\noindent J. Yadav and G. Shanker\newline
Department of Mathematics and Statistics\newline
Central University of Punjab\newline
Bathinda, Punjab-151401, India.\newline
Email: sultaniya1402@gmail.com; gauree.shanker@cup.edu.in\newline\\
				%
				
				%
				%

			\end{document}